\documentclass{cmslatex}
          
\usepackage[paperwidth=7in, paperheight=10in, margin=.8in]{geometry}

\usepackage{amsmath,amsfonts,amssymb}
\usepackage{enumitem}
\usepackage{hyperref}
\usepackage{lineno}
\usepackage{color}
\usepackage{graphicx,epsfig}
\usepackage[percent]{overpic}
\usepackage{tcolorbox}
\usepackage{tikz-cd}
\usepackage{mathrsfs}
\usepackage[ansinew]{inputenc}
\usepackage{float}
\usepackage{afterpage}
\usepackage{color}

\usepackage{subfigure}

\newcommand{\bX}{\textbf{X}}
\newcommand{\bV}{\textbf{V}}
\newcommand{\bF}{\textbf{F}}
\newcommand{\ks}{k^*}
\newcommand{\kss}{k^{**}}
\renewcommand{\S}{\mathcal S}

\renewcommand{\theequation}{\arabic{section}.\arabic{equation}}

\begin{document}

\title{Preventing Congestion in Crowd Dynamics Caused by Reversing Flow}

\author{
Giuseppe G. Amaro\thanks{
        GAe Engineering S.r.l., Via Assietta 17, 10128 Turin, Italy.}
\and Emiliano Cristiani\thanks{
	Istituto per le Applicazioni del Calcolo ``M.\ Picone'', 
	Consiglio Nazionale delle Ricerche,
	Via dei Taurini, 19 --
	00185 Rome, Italy, 
	e.cristiani@iac.cnr.it (corresponding author)}
\and Marta Menci\thanks{
	Istituto per le Applicazioni del Calcolo ``M.\ Picone'', 
	Consiglio Nazionale delle Ricerche, Rome, Italy, 
	m.menci@iac.cnr.it}
	}
	
\pagestyle{myheadings} 
\markboth{Preventing Congestion in Crowd Dynamics Caused by Reversing Flow}{G.G. Amaro, E. Cristiani, M. Menci} 
\maketitle



\begin{abstract}
\phantomsection\addcontentsline{toc}{section}{\numberline{}Abstract}

In this paper we devise a microscopic (agent-based) mathematical model for reproducing crowd behaviour in a specific scenario: a number of pedestrians, consisting of numerous social groups, flow along a corridor until a gate located at the end of the corridor closes. 
People are not informed about the closure of the gate and perceive the blockage observing dynamically the local crowd conditions. 
Once people become aware of the new conditions, they stop and then decide either to stay, waiting for reopening, or to go back and leave the corridor for ever. 
People going back hit against newly incoming people who are not yet aware of the blockage or have already decided to stay. This creates a dangerous counter-flow which can easily lead to accidents.
We run several numerical simulations varying parameters which control the crowd behaviour, in order to understand the factors which have the greatest impact on the dynamics. We conclude with some useful suggestions directed to the organizers of mass events. 

\end{abstract}

\begin{keywords}
Crowd modelling, crowd control, social force model, counter-flow, social groups
\end{keywords}

\begin{AMS}
76A30
\end{AMS}

\section{Introduction}\label{sec:intro}
This paper proposes a microscopic (agent-based) model for describing the dynamics of pedestrians in a very special situation characterized by a moving crowd whose members gradually reverse their direction of motion, thus colliding with the newly incoming people still moving in the original direction.
This situation creates a nonstandard self-regulating counter-flow with a rather complex dynamics, which can be investigated by suitably tuning the parameters of the model. 
The final goal is to predict (and prevent) the formation of dangerous congestion in real mass events.

\medskip 

\emph{Relevant literature.} 
Modelling the behavior of people in a crowd is a difficult task, since it requires to identify the most important behavioral rules which greatly vary from person to person. For that reason, the study of crowds can be regarded as a multidisciplinary area, which have attracted since many years the interest of mathematicians, physicists, engineers, and psychologists.

Crowd modelling has a long-standing tradition, starting from the pioneering papers by 
Hirai and Tarui \cite{hirai1975}, 
Okazaki \cite{okazaki1979TAIJa}, 
and Henderson \cite{henderson1974} in the '70s.
Since then, all types of models were proposed, spanning from microscale to macroscale, including multiscale ones, both differential and nondifferential (e.g., cellular automata).
Models can be first-order (i.e.\ velocity based) or second-order (i.e.\ acceleration based), with local or nonlocal interactions, with metric or topological interactions, with or without contact-avoidance features. The presence of social groups can also be taken into account.
A number of review papers  
\cite{aghamohammadi2020, 
	bellomo2011,
	dong2020, 
	duives2013,
	eftimie2018,
	haghani2020, 
	martinez2017,
	papadimitriou2009}, 
and books \cite{cristiani2014book, rosini2013book, kachroo2008book, maurybook2019} are now available,  
we refer the interested reader to these references for an introduction to the field. 
It is also useful to mention that models for pedestrians often stem from those developed in the context of vehicular traffic \cite{helbing2001, rosini2013book}.
Moreover, there is a strict connection between pedestrian modeling and control theory, including mean-field games, see, e.g., \cite{albi2020, cristiani2021pp, cristiani2014book} and reference therein.

In this paper we deal specifically with the \emph{counter-flow} dynamics: such a dynamics occur when two groups of pedestrians move in opposite directions, therefore each group has to find a way to pass through the other. 
The importance of accurate modeling and simulation of counter-flow dynamics is supported by evidence from crowds disaster analysis \cite{helbing2012crowd}.  
Over the last years, the major accidents often occurred for overcrowding and counter-flow phenomena taking place inside the area of the events, or in the proximity of the entrances and exit points. 
For all these reasons, the literature on counter-flow is quite reach; see, among others, \cite{helbing1995social, hoogendoorn2003simulation, heliovaara2012counterflow} in the context of microscopic differential models, \cite{cristiani2011MMS} for a multiscale differential models, and \cite{weng2006cellular, nowak2012quantitative} in the context of cellular automata. 
Moreover, it is now well established that counter-flow dynamics lead to the so-called \emph{lane formation}: in order to avoid collision, pedestrians arrange in alternate lanes ($\leftrightarrows$) having the same walking direction.
The behavior displayed by numerical simulations is in good agreement with observations of real people in both artificial and natural environments \cite{Kretz2006experimental, hoogendoorn2003extracting, helbing2001self, murakami2021mutual}.

Finally, let us discuss the impact of the \emph{social groups} on the dynamics of crowds: despite most of the models assume that each pedestrian moves in a crowd on its own, real crowds are typically formed by small subgroups, such as friends or families.
The impact of social groups on crowd dynamics has been explored since the '70s in a number of papers \cite{aveni1977not, moussaid2010walking, von2017empirical, singh2009modelling}, resulting in the fact that the presence of groups is not negligible at all. 
In particular, in the context of bi-directional counter-flow dynamics, theoretical and experimental observations suggest that the presence of groups slows down the formation of lanes, which are more fragmented \cite{zanlungo2020effect, crociani2017micro}.

\medskip

\emph{Paper contribution.}
In this paper, we consider a specific scenario related to the counter-flow dynamics: a crowd (with social groups) in a corridor, initially moving in one direction towards an open gate, at some moment is no longer able to proceed because of the closure of the gate. After that, people have to decide whether to stop \& wait for a possible re-opening, or to move back.
This is the case, e.g., of a inflow of people towards an area dedicated to a mass event, which is interrupted by the organizers when the area capacity has been reached.

From the modelling point of view, the main novelty is that the decision to stay or to move back is taken \emph{dynamically} by each group, on the basis of the behaviour of the surrounding groups. More precisely, we assume that the decision is taken whenever \emph{the group leader is no longer able to move forward for a certain time}. 
This can happen either because it has reached the gate or it has reached the other people queuing in front of the gate, or it is hit by the people going back and blocking the way. 
The overall dynamics is also complicated by the fact that people who decide to stay do not want to be overtaken by the other people approaching the gate, since they do not want to lose their priority in the queue. On the other hand, staying people want to facilitate the passage of reversing people, because the latest are leaving free space for the former.

From the mathematical point of view, we propose a microscopic differential model inspired by the well-known Helbing's Social Force Model \cite{helbing1995social}, based on a large system of ODEs. 
We introduce several variants with respect to the original model: the two most important of them are that a) we consider a first-order (velocity-based) model and b) we consider topological interactions, taking into account the first neighbour only. The last choice greatly speed up the numerical code. 

The final goal of this research, which is especially dedicated to practitioners who are involved in the organization of real mass events, is the prevention of critical situations which could arise in the reversing-flow scenario under consideration. 
Moreover, since digital twins are commonly used to support the safety plan development, highlighting critical aspects that need to be solved before and during the event (cf.\ \cite{scozzari2018modeling}), we also propose a crowd control strategy based on the optimal placing of signals/stewards, aiming at informing people in due time about the gate status (open or close). 

\medskip

\emph{Paper organization.}
The paper is organized as follows. 
In Section \ref{sec:model} we present the model for the crowd dynamics. 
In Section \ref{sec:casestudy} we present our case study: we describe the geometry and the mechanism undergoing the behavioural choices of pedestrians. 
In Section \ref{sec:tests} we discuss the results obtained from numerical simulations. 
We end the paper with some conclusions and future perspectives.

\section{The model}\label{sec:model}
\subsection{General principles}
All force models share the same structure of the Newtonian dynamics, namely
\begin{equation}\label{structure}
\left\{
\begin{array}{l}
\dot{\bX}_k(t)= \bV_k(t) \\ [2mm]
\dot{\bV}_k(t)= \bF_k(t,\bX,\bV)
\end{array},
\qquad k=1,\ldots,N
\right.
\end{equation}
where $N$ is the total number of agents, $\bX_k(t),\bV_k(t)\in \mathbb{R}^{2}$, denote the position and velocity of agent $k$ at time $t$, respectively, and $\bX=\left(\bX_1,\ldots,\bX_N\right)$, $\bV=\left(\bV_1,\ldots,\bV_N\right)$. 
The function $\mathbf{F}_k$, the so-called \emph{social force}, models the total force exerted on agent $k$, and gathers all the physical, psychological and behavioral aspects of pedestrian dynamics.
The social force is not a real force, but rather an empirical mathematical tool which translates in formulas all these aspects. In its minimal form, it takes into account the following three contributions:
\begin{itemize}
	\item[i)] An individual desired velocity term: the velocity that a single pedestrian would keep if it was alone in the domain.
	\item[ii)] A repulsion term: pedestrians tend to maintain a certain distance with respect to other members of the crowd, obstacles and walls present in the environment, in order to avoid collisions. 
	\item[iii)] An attraction term: pedestrians who are not moving on their own tends to stay close to the members of their social group (friends, family members). 
\end{itemize}

Social force models can be enhanced adding random fluctuations due to unpredictable behavioral variations or further small-scale interactions. In this regards, it can be useful to note that the (generally undesired) numerical instabilities often play the same role. 

\medskip

\emph{First-order models.} 
In pedestrian dynamics, much more than in vehicular dynamics, accelerations are almost instantaneous (at least if compared with a reference typical time scale) and then inertia-based effects are negligible. 
Therefore, first-order models of the form
\begin{equation}\label{structure-firstorder}
\dot\bX_k(t)=\bV_k(t,\bX),\qquad k=1,\ldots,N
\end{equation}
are also suitable, see, e.g., \cite{cristiani2011MMS}.  
We think that such a models are easier to calibrate and are computationally less expensive, for these reasons in this paper we will adopt a model of this kind.
This means that all the aspects of the dynamics which were encapsulated in the force $\bF$ are now inserted directly in the velocity vector $\bV$.

\medskip 

\emph{Social groups.}
In order to account for the presence of social groups in the crowd, such as families or friends, we assume that each pedestrian is part of a group. 
Each group has at least two members, we do not consider the presence of lonely people.
Moreover, each group has a leader, which never changes in the time frame of the simulation. 
The leader of the group takes decisions about the common target of the whole group. 
Groups tend to stay together, but they can temporarily break up and then reunite.
The leader does not necessarily walk in front of the group because all group members know the destination (decided by the leader) and are able to reach it independently.

\medskip

\emph{Groups status.}
We assume that at any given time, each group (identified with its leader) has a unique \emph{status} which corresponds to its target and, more in general, to its behaviour. The four possible statutes will be detailed later on in Sect.\ \ref{sec:behavior}.

\medskip

\emph{Topological interactions.}
We consider topological, rather than metric, interactions, meaning that each agent $k$ interacts with a fixed number of agents at the same time, regardless of their distance from the agent $k$. 
More precisely, we assume that each leader interacts with the first neighbour outside its social group only, while followers (i.e.\ not leaders) interact with the first neighbour inside their social group, the first neighbour outside, and their leader.
As it happens in molecular dynamics, the fact that a particle interacts with a few other particles at a time does not mean that the interactions are limited to them. The first neighbour changes continuously and after few time iterations all the agents reasonably close to each agent have interacted with it. 
While standard social force model assumes contemporary interactions with neighbours within a certain threshold distance, we prefer to consider fewer interactions at a time, recovering the same results, in average, over a longer, but still short, time period.
This choice stabilizes the dynamics and it is convenient from the computational point of view.

\subsection{Mathematical details}
As we said, we consider a first-order social force model with social groups and topological interactions. Agents have no dimension and we do not consider contact-avoidance features.

Let us start with the dynamics of \emph{group leaders}. 
We set
\begin{equation}\label{leaderdynamics}
\bV_k(\bX;\S)=
\bV^d(\S_k)+
\bV^{R}(\bX_k,\bX_{\kss};\S_k,\S_{\kss})+
\bV^o(\bX_k)
\end{equation}
where
$\bV^d$ is the desired velocity, which only depends on the status $\S_k$ of the group the agent $k$ belongs to;
$\bV^{R}$ accounts for the repulsion from strangers, and depends on the positions $\bX_k$ of the agent $k$ itself, the position $\bX_{\kss}$ of the agent $\kss$, defined as the nearest neighbour of $k$ \emph{outside} its social group, and the statuses of the two agents; 
$\bV^o$ is the repulsion from obstacles.
Clearly we have $\S:=(\S_1,\ldots,\S_N)$. 

The dynamics of \emph{followers} is instead 
\begin{equation}\label{followerdynamics}
\bV_k(\bX;\S)=
\bV^d(\S_k)+
\bV^{r}(\bX_k,\bX_{\ks})+
\bV^{R}(\bX_k,\bX_{\kss};\S_k,\S_{\kss})+
\bV^a(\bX_k,\bX_{k^L})+
\bV^o(\bX_k)
\end{equation}
where 
$\bV^{r}$ accounts for the repulsion from group members, and depends on the positions $\bX_k$ of the agent $k$ itself, the position $\bX_{\ks}$ of the agent $\ks$, defined as the nearest neighbour of $k$ \emph{inside} its social group;
$\bV^a$ accounts for the attraction towards the group's leader, whose index is denoted by $k^L$.

Note that group leaders are not attracted by group mates, the cohesion of the group being totally left to followers.
Moreover, leaders are not repulsed by group mates. This avoids artifacts in the dynamics and self-propulsion of the group.

The reason why we distinguish internal and external repulsion is that pedestrians tend to stay close to group mates, while keeping a larger distance from the strangers. 

\medskip

The repulsion is defined in order to be inversely proportional to the distance between the agents,
\begin{equation} 
\bV^{r}= -C^{r}\frac{\bX_{\ks} - \bX_k}{\left\| \bX_{\ks} - \bX_k \right\|^2},
\qquad
\bV^{R}= -C^{R}(\S_k,\S_{\kss})\frac{\bX_{\kss} - \bX_k}{\left\| \bX_{\kss} - \bX_k \right\|^2},
\end{equation}
where the parameter $C^{r}>0$ is constant, whereas $C^{R}(\S_k,\S_{\kss})>0$ depends on the behavioral statuses of the two interacting agents, see Sect.\ \ref{sec:behavior}.
To avoid numerical issues in the case two agents get temporarily too close to each other, a threshold is imposed to the value $1/\|\cdot\|^2$. 

The attraction follower-leader, instead, is proportional to the distance between the two agents
\begin{equation} 
\textbf{V}^a= C^a
(\bX_{k^L} - \bX_k),
\end{equation}
where the parameter $C^{a}>0$ is constant.


\section{Description of the case study} \label{sec:casestudy}

\subsection{Geometry}
\begin{figure}[h!]
	\centering
	\includegraphics[width=12cm, height=2.5cm]{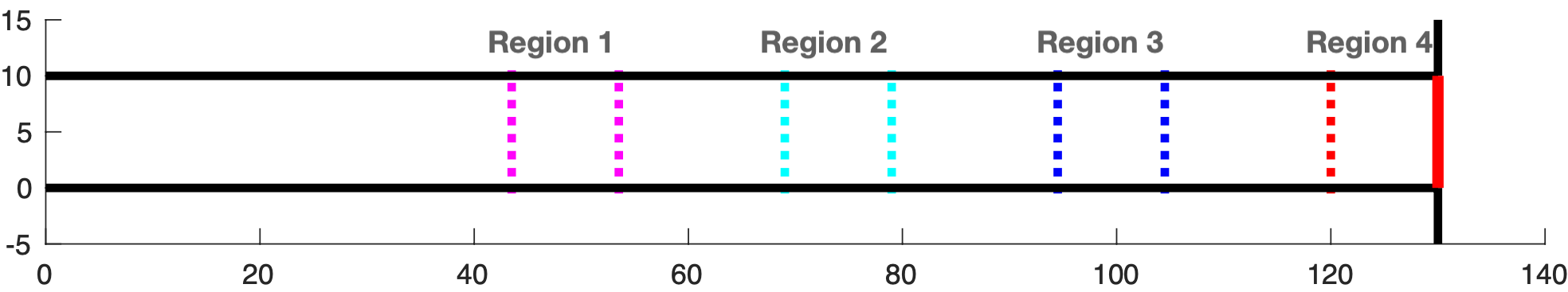}
	\caption{Geometry of the simulated domain. The gate is placed on the far right. Pedestrian density evolution is monitored within the four square regions placed along the corridor.}
	\label{fig:geometry}
\end{figure}

We model the access road to the area of a mass event as a two-dimensional corridor of length $L$ and height $H$.
People initially move along the corridor from left to right and enter the area of the event through a gate, placed at the end of the corridor, see Fig.\ \ref{fig:geometry}. At a certain time ($t=0$ in the simulation) the gate closes for safety reasons so that people are forced to change target (staying or going back). 
In the following numerical simulations we focus on the behavior of the crowd after the gate closure, which involves a risk of high densities. 
In particular, we measure the pedestrian density along the corridor by means of four square regions of side $R$, equally spaced along the corridor. Regions are numbered from left to right, see Fig.\ \ref{fig:geometry}.

\subsection{Behavior}\label{sec:behavior}
In order to complete the description of the model we need to specify how the groups take decisions. In particular, how and when they decide to move towards the gate, stop or go back.

We select 4 statuses $\S_k\in\{1,2,3,4\}$ at the group level (i.e.\ all components of a group share the same status, decided by the group leader), defined as follows:
\begin{itemize}
	\item $\S_k=1$: move rightward towards the gate;
	\item $\S_k=2$: doubt phase, decision in progress;
	\item $\S_k=3$: decision taken, go back moving leftward;
	\item $\S_k=4$: decision taken, queuing in the corridor.
\end{itemize}
The group status affects the desired velocity and the interactions with the others.

Status' changes happen necessarily in this order: 
$1 \rightarrow 2 \rightarrow \{3,4\}$, and once the group is in status 3 or 4 it will not change further.

The change 1 $\rightarrow$ 2 occurs if    
\begin{equation}
t>\delta t \quad \text{ and } \quad X^1_k(t)-X^1_k(t-\delta t)\leq\overline{\delta\ell},
\end{equation}
where $\delta t$, $\overline{\delta \ell}>0$ are two additional parameters and $X^1_k$ is the horizontal component of the position of the agent $k$. 
This means that the group leader was not able to move forward for more than $\overline{\delta\ell}$ in a period of time of length $\delta t$. In practice, $\delta t$ represents a degree of willingness to get to the gate.

The status 2 has a fixed duration, set to $D$. 
After $D$ time units, the change 2 $\rightarrow$ $\{3,4\}$ occurs randomly in such a way that, in average, $p$\% of groups fall in status 3 and $(1-p)$\% in status 4.

\section{Numerical tests}\label{sec:tests}

\subsection{Parameters}
Let us begin with setting the values of the parameters. Parameters are primarily chosen in order to get realistic (observed) behaviour in situations of equilibrium. 
In particular we match observed crowds in the initial phase, where all people are in status 1 and are normally walking rightward, and in the final phase, when only people with status 4 are present and waiting. 
In these two phases the forces which rule the interpersonal distances are in equilibrium, and the relative positions of persons do not vary (apart for small negligible oscillations).

In the numerical tests we discuss the effect of five parameters, namely $N$, $\delta t$, $p$, $D$, $H$. 
We also investigate the effect of the desired velocity of pedestrians with status 4.

Equation \eqref{sec:model} is numerically solved by means of the explicit Euler scheme, with time step $\Delta t$.

\medskip

\emph{Fixed parameters.}
$\Delta t=0.01$ s,  
$\overline{\delta\ell}=1.5$ m, 
$L=130$ m, 
$R=10$ m,
$\bV^d(1)=(1,0)$ m/s,
$\bV^d(2)=(0,0)$ m/s,
$\bV^d(3)=(-1.2,0)$ m/s.
The 16 possible values of $C^R$ are summarized in the following matrix, where the entry $ij$ is the force exerted by and agent in status $j$ over an agent in status $i$,

$$
\left[
\begin{array}{cccc}
2.0  & 2.5  & 1.0  & 2.5  \\
2.0  & 2.0  & 2.0  & 2.0  \\
0.75 & 0.75 & 0.75 & 0.75  \\
2.0  & 2.0  & 0.75 & 4.5  
\end{array}
\right]
$$
Some comments on the repulsion forces are in order:
\begin{enumerate}
	\item People in status 3 (coming back) are repelled from all the others with little intensity. In fact, they accept to be close to other people, even to strangers, because the contact is supposed to be of short duration. 
	\item People in status 4 (staying) are strongly repelled by people of the same status because they want to reach a large comfort distance. 
	\item People in status 4 (staying) are also repelled by people in status 1 (moving rightward) because they want to keep their priority in the queue and keep the newly coming people behind them. In other words, they do not leave room for others to move rightward. Conversely, they are little repelled by people in status 3 (coming back) because waiting people take advantage from the departure of them. People going back leave free space and allow waiting people to get even closer to the gate.
	All these features have an important consequence: staying people are compressed towards the gate and the closer they are to the gate, the more compressed they are.
\end{enumerate}

\emph{Variable parameters.}
$N$ = 400 or 800 or 1200,  
$\delta t$ = 7 or 20 s, 
$p$ = 25 or 75, 
$D$ = 10 or 40 s, 
$\bV^d(4)$=$(0.5,0)$ or (0,0) m/s,
$H$ = 7 or 10 or 13 m.

\medskip

In all simulations the initial density of pedestrians is 0.8 p/m$^2$ (the larger $N$ the larger the area occupied).
The initial positions are randomly chosen.
The number of members of the groups is uniformly random in the set $\{2,3,4,5,6\}$.

\subsection{Numerical results}

Fig.\ \ref{fig:screenshot} shows some screenshots of a reference simulation with 1200 agents: 
\begin{figure}[t!]
	\centering
	{\footnotesize (a)}\subfigure{\includegraphics[width=13.3cm, height=1.9cm]{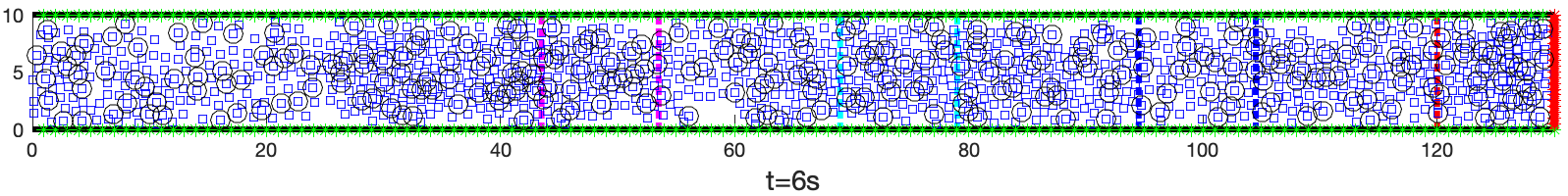}}\\
	{\footnotesize (b)}\subfigure{\includegraphics[width=13.3cm, height=1.9cm]{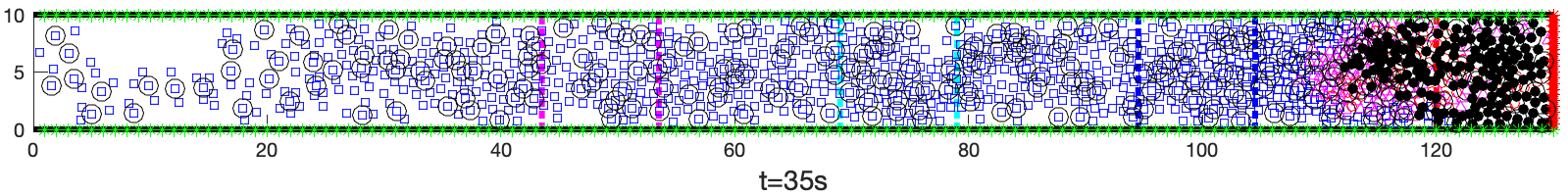}}\\
	{\footnotesize (c)}\subfigure{\includegraphics[width=13.3cm, height=1.9cm]{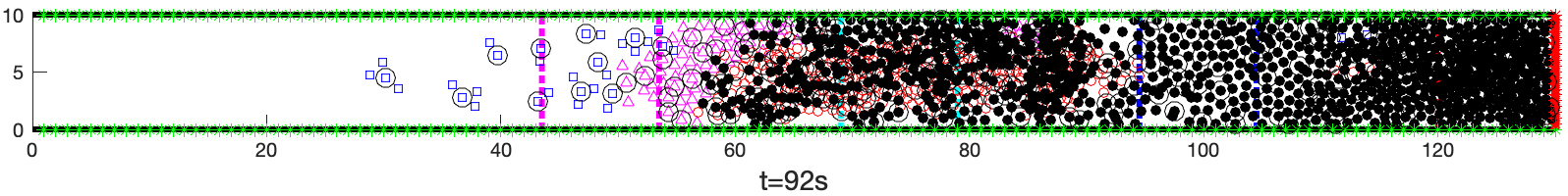}}\\
	{\footnotesize (d)}\subfigure{\includegraphics[width=13.3cm, height=1.9cm]{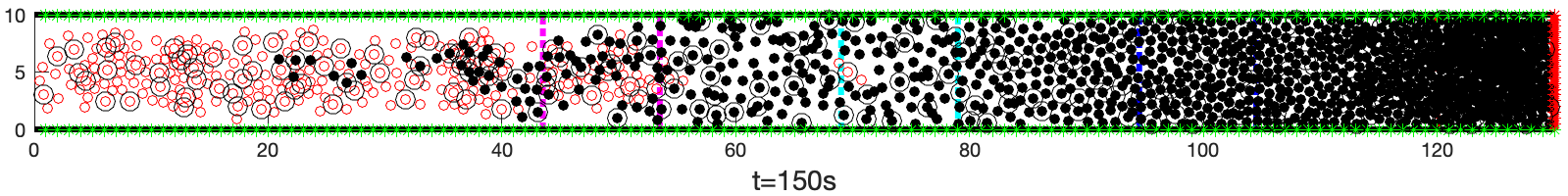}}\\
	{\footnotesize (e)}\subfigure{\includegraphics[width=13.3cm, height=1.9cm]{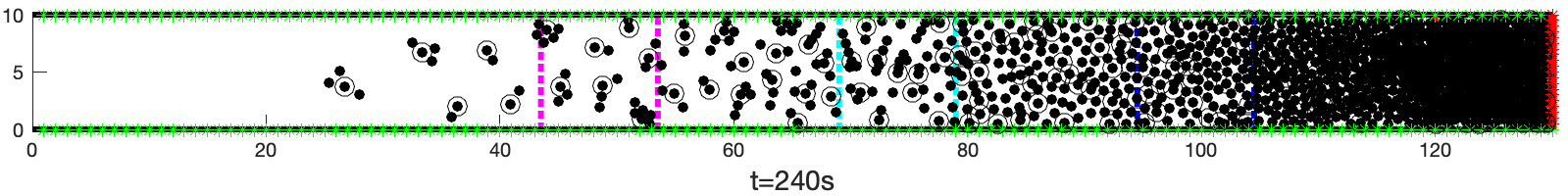}}
	\caption{Five screenshots of a simulation with $N=1200$ pedestrians, with $\delta t$ = 7 s, $p$=25, $D=10$ s, $\bV^d(4)$=$(0.5,0)$ m/s, $H$ = 10 m.
		Status 1 is blue, status 2 is magenta, status 3 is red, and status 4 is black. Circled agents are the leaders of their groups.}
	\label{fig:screenshot}
\end{figure}
(a) at the beginning all agents are walking rightward. As we said before, initial density, repulsion forces and desired velocity are compatible with each other, in the sense that people can move forward at constant velocity, nobody slows down due to excessive proximity to the neighbor and no queue is formed.
As soon as the first pedestrians reach the closed gate the situation changes: people in proximity of the gate stop and people behind start slowing down. (b) A queue begins to form and gradually still people take decision about the new status to get.
At this point all four statuses are present at the same time and the dynamics becomes quite complex. 
Some people in status 1 (moving rightward, blue) are able to overcome people in status 4 (staying, black) before changing status themselves.
(c) Lanes are formed as in more standard counter flow dynamics (cf.\ Sect. \ref{sec:intro}), although they are quite perturbed by the presence of the social groups. 
Social groups sometimes disperse for a while but then reunite.
(d) Gradually people in status 3 (coming back, red) succeed in passing through the staying crowd and leave the corridor. 
While they move leftward, people in status 4 (staying, black) can move a bit backward to let them free space.
(e) Finally, only people in status 4 (staying, black) are present in the corridor. They reach an equilibrium characterized by a noncostant density, due to the natural tendency to stay close to the gate.

\medskip

\begin{figure}[t!]
	\centering 
	\subfigure[][$N_p=400$]{\includegraphics[width=4.2cm, height=4.2cm]{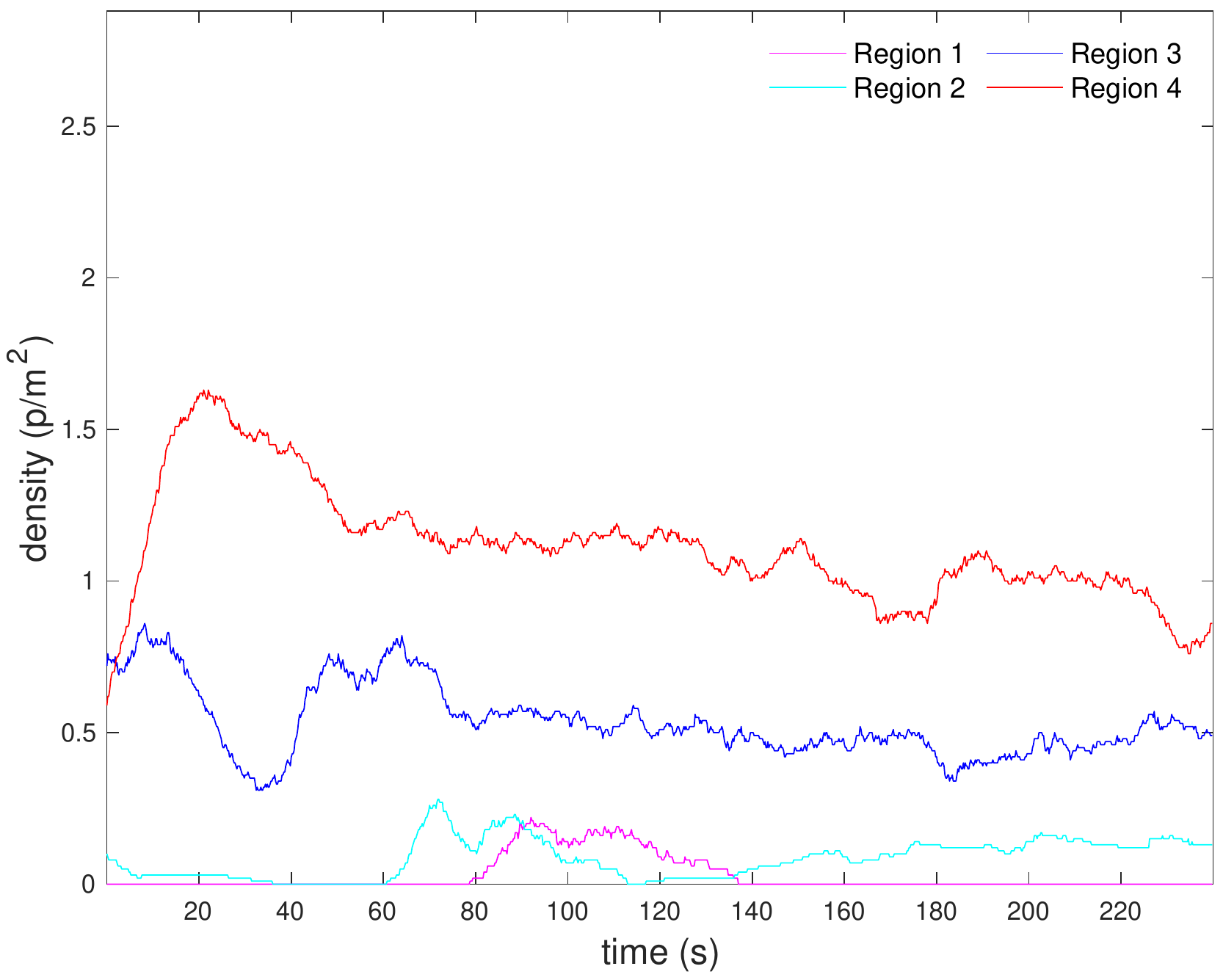}}
	\hspace{0 mm}
	\subfigure[][$N_p=800$]{\includegraphics[width=4.2cm, height=4.2cm]{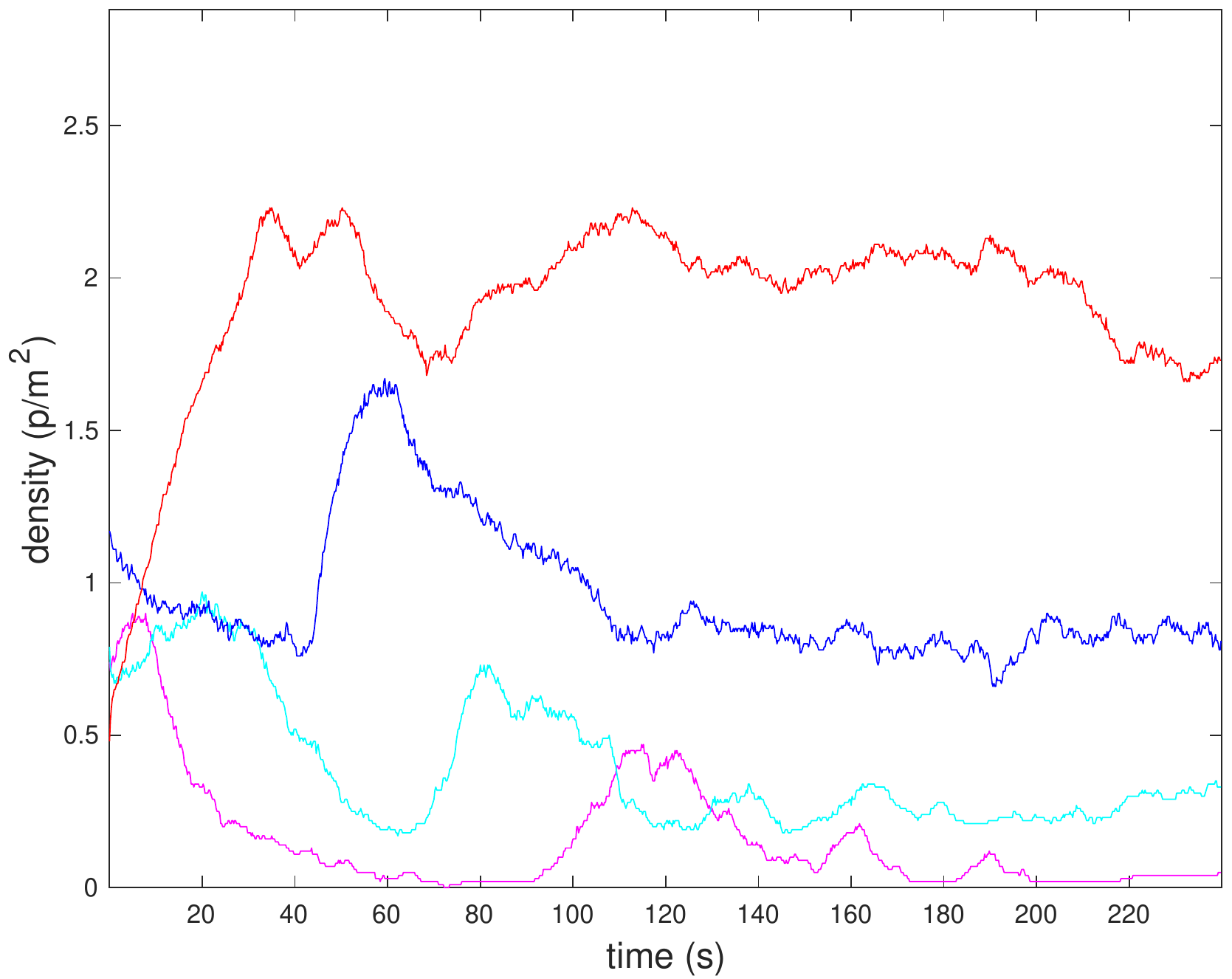}}
	\hspace{0 mm}
	\subfigure[][$N_p=1200$]{\includegraphics[width=4.2cm, height=4.2cm]{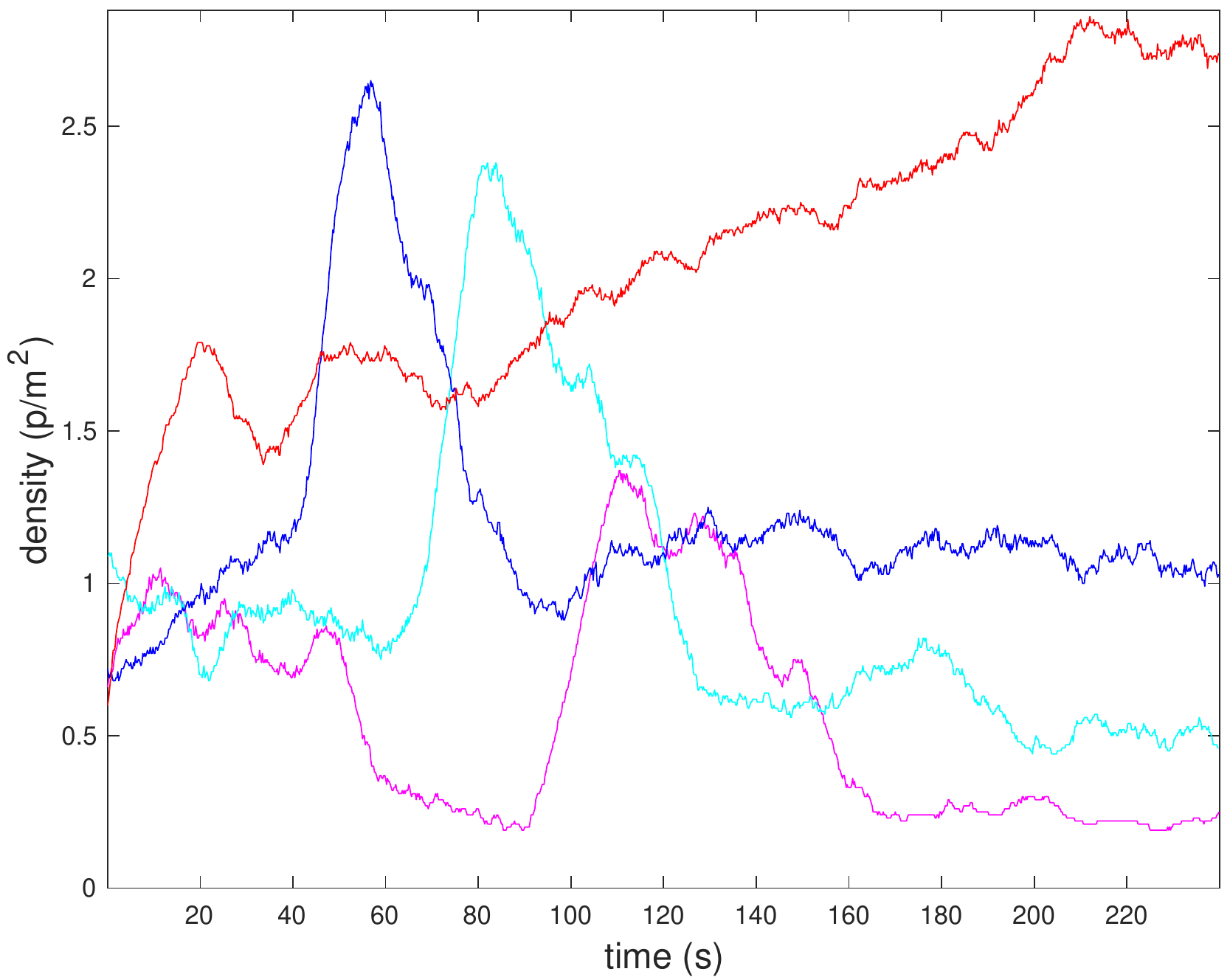}}
	\caption{Evolution in the four regions of interest, for $N_p=400, 800, 1200$ and other parameters as in Fig.\ \ref{fig:screenshot}.}
	\label{fig:densities_new}
\end{figure}

Fig.\ \ref{fig:densities_new} shows instead the evolution of the density in the four regions of interests, varying $N$. 
The comparison of the plots brings to light some facts: 
\begin{enumerate}
	\item The most crowded region is the fourth one, the region closest to the gate.
	\item In all cases it is well visible a peak of the density in regions 3-2-1 which decreases and shifts forward in time as the region moves away from the gate (see, e.g., Fig.\ \ref{fig:densities_new}(c)). This moving peak is due to the reversing people encountering the newly arriving people.
	\item As $N$ increases, average densities increase in all regions but less than one could expect. We think that the reason for this stability is that \emph{the crowd is able to self-regulate}: the more people in the corridor, the earlier (i.e.\ before in both space and time) people get to make a decision about what to do, and, if they decide to leave, they do so before remaining trapped by the others. We observe a sort of a compensation phenomenon that consists in the fact that the more people there are, the sooner they leave, overall leading to small differences in the density evolution.
\end{enumerate}

\begin{figure}[t!]
	\centering
	\subfigure[][$H=7$ m]{\includegraphics[width=4.2cm, height=4.2cm]{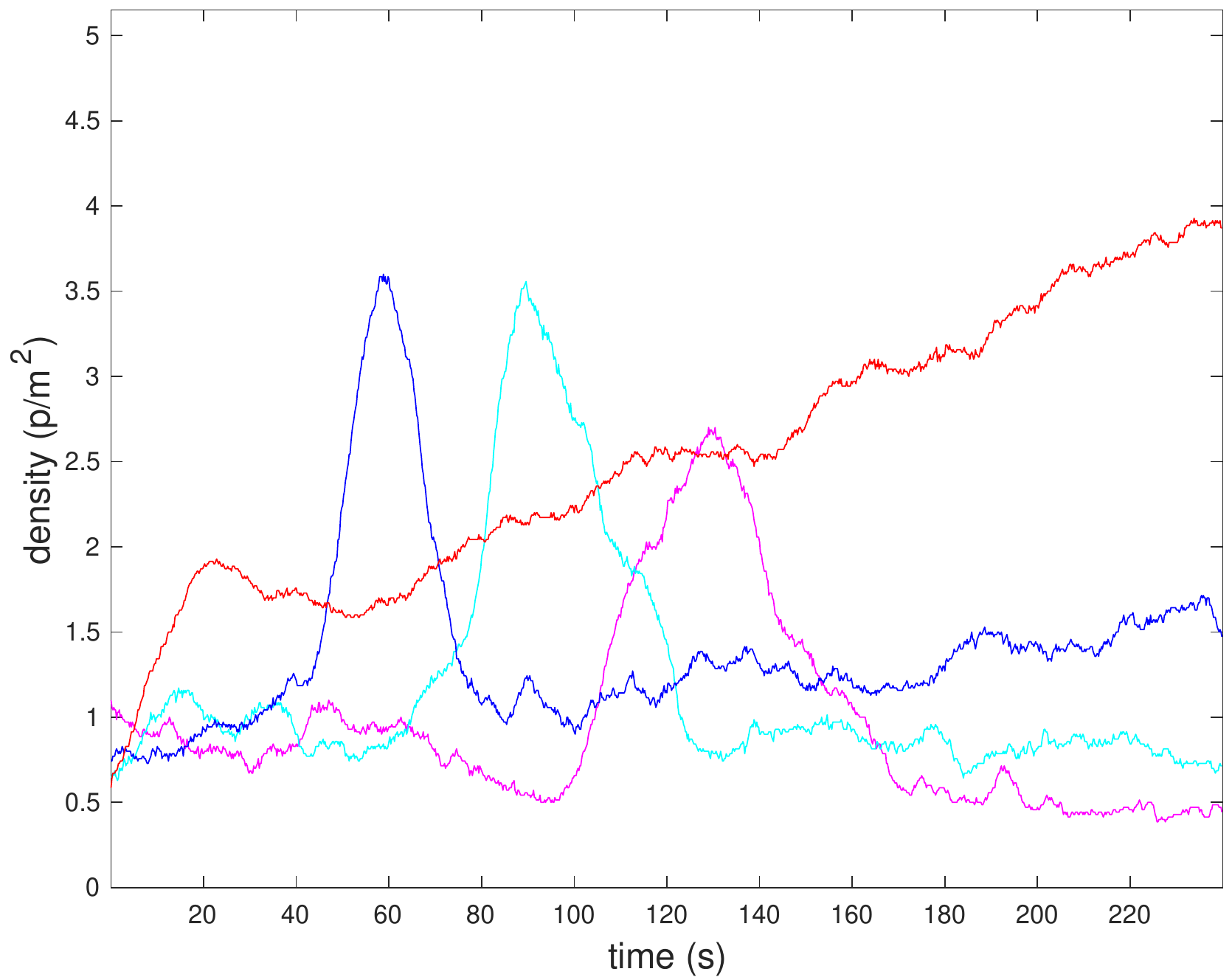}}
	\hspace{0 mm}
	\subfigure[][$H=13$ m]{\includegraphics[width=4.2cm, height=4.2cm]{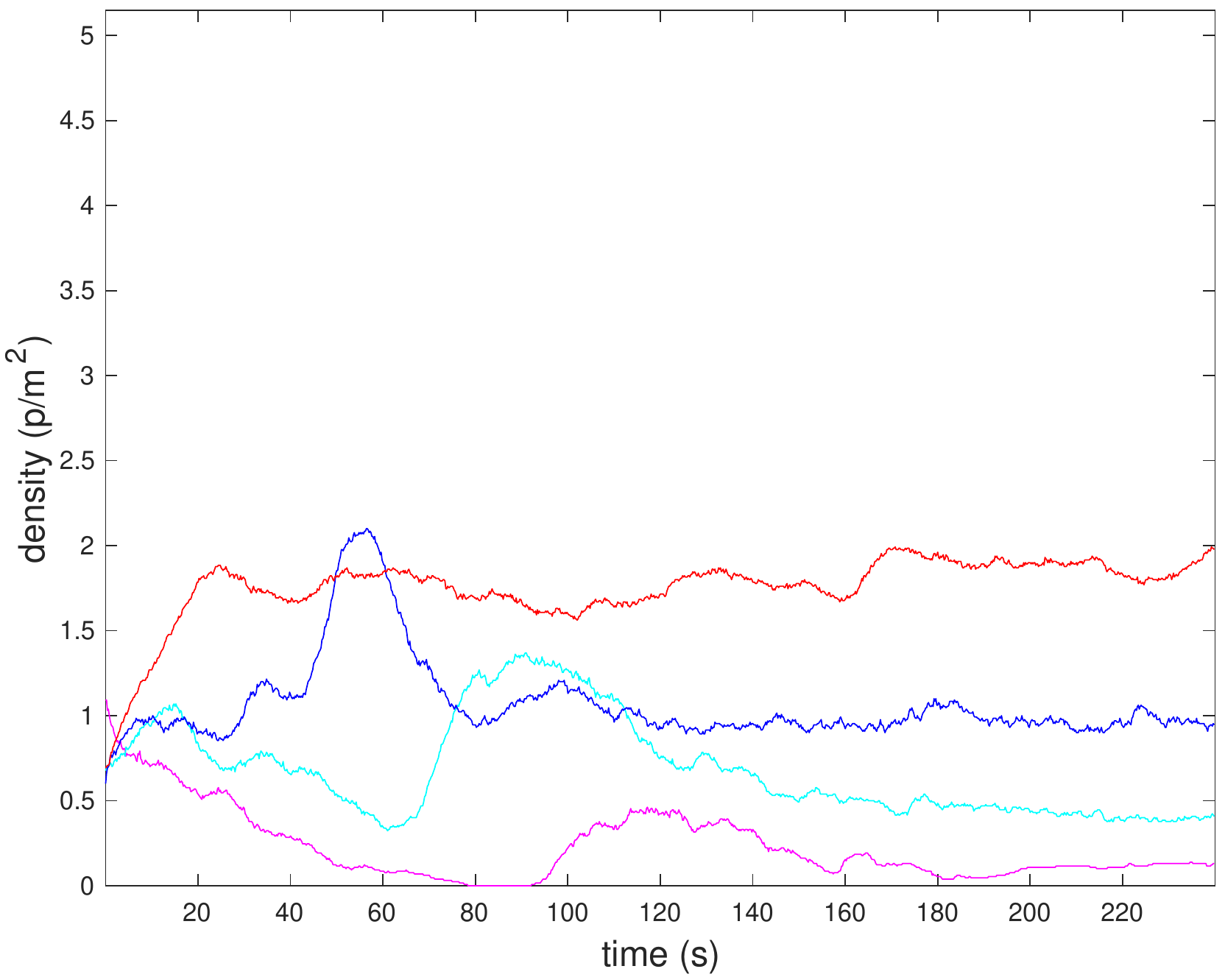}}
	\hspace{0 mm}
	\subfigure[][$D=40$ s]{\includegraphics[width=4.2cm, height=4.2cm]{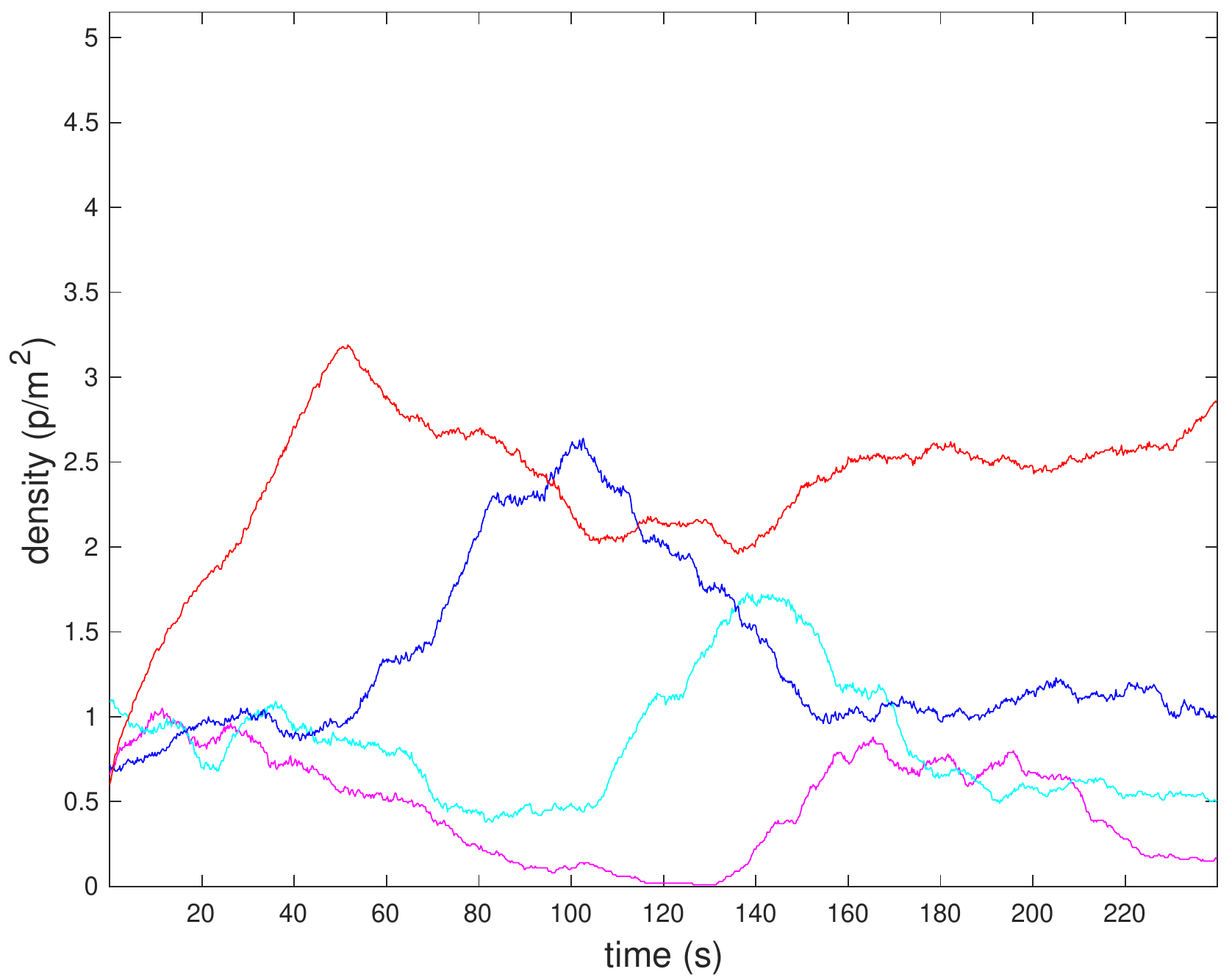}}\\
	\subfigure[][$\delta t=20$ s]{\includegraphics[width=4.2cm, height=4.2cm]{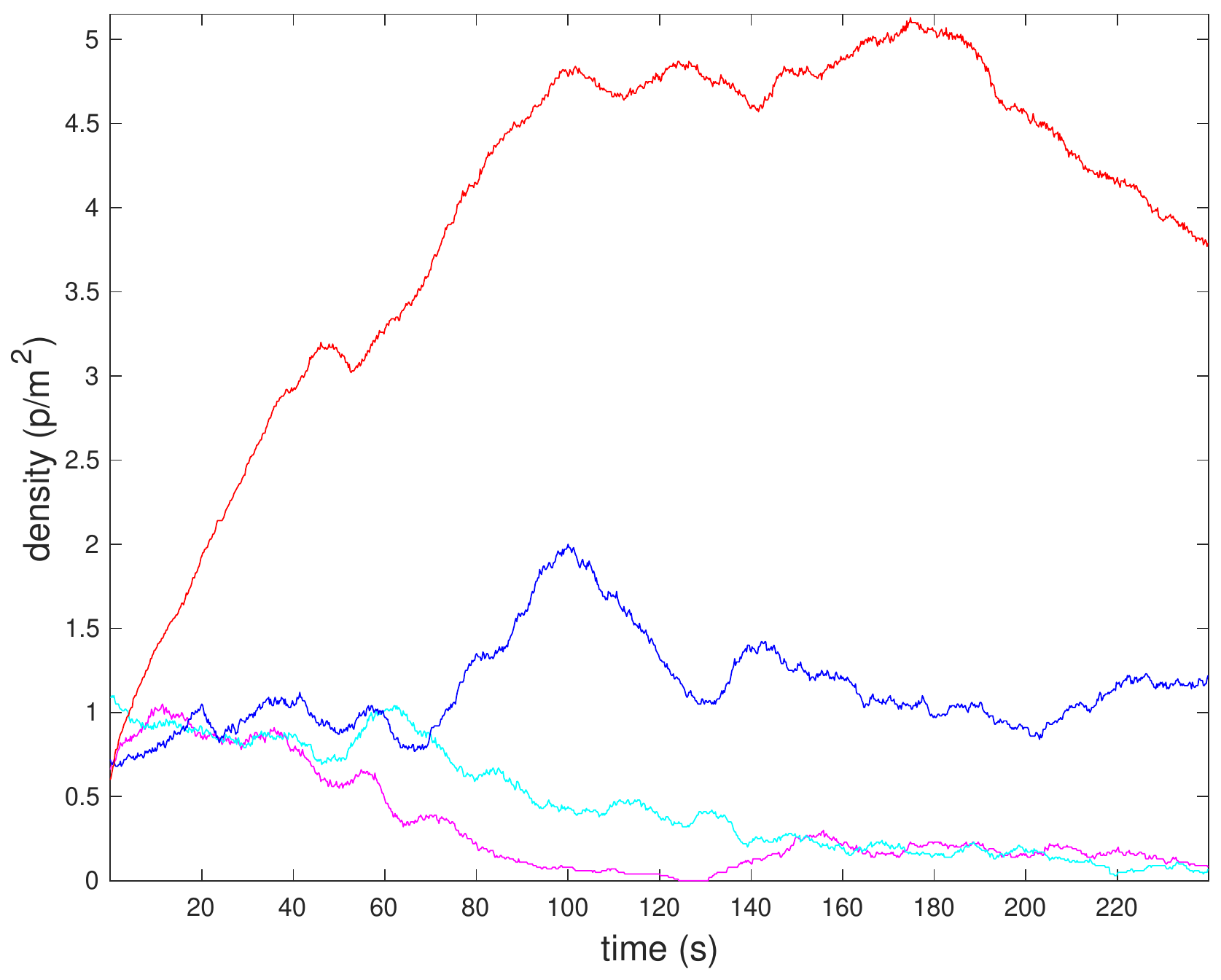}}
	\hspace{0 mm}
	\subfigure[][$p=75$]{\includegraphics[width=4.2cm, height=4.2cm]{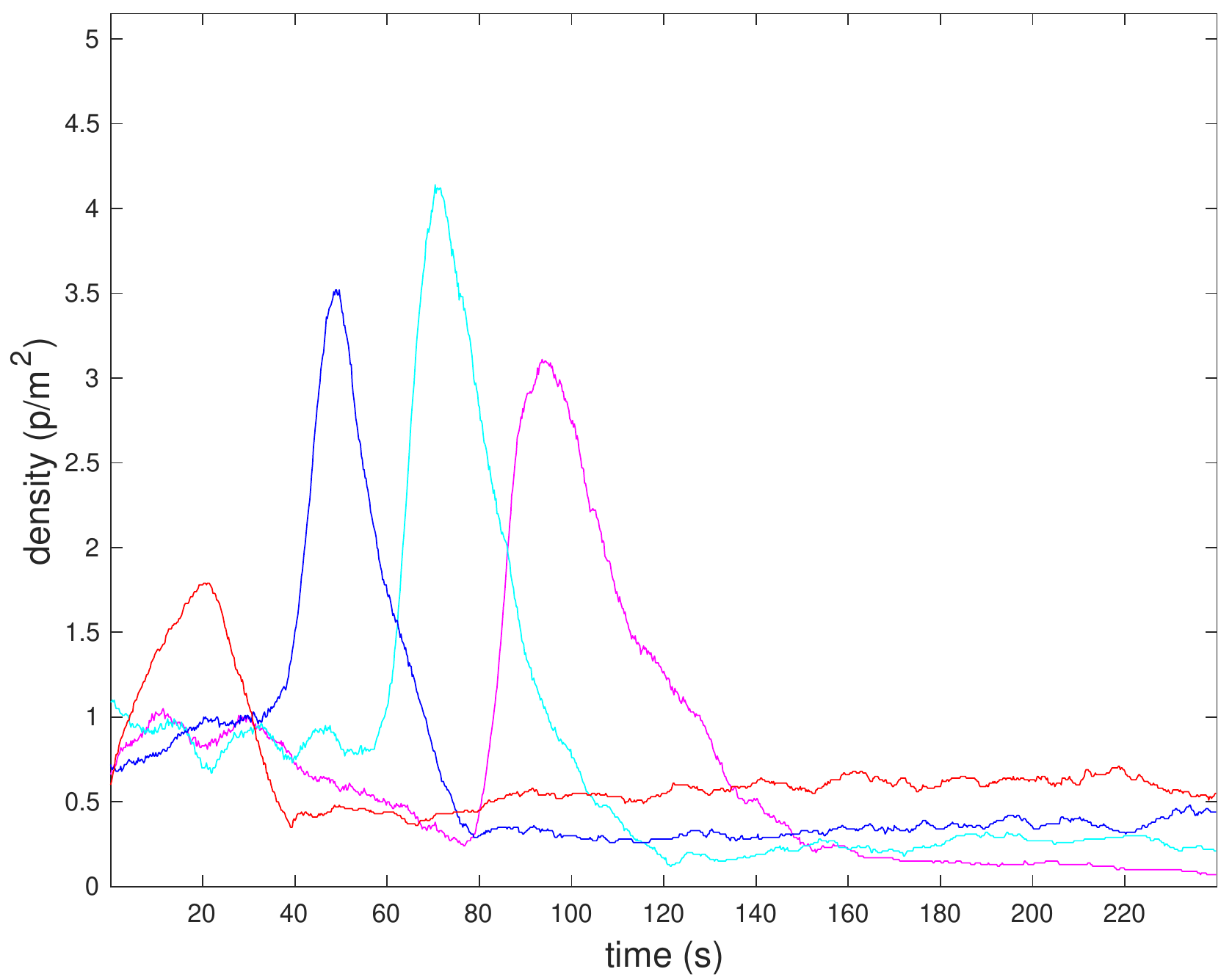}}
	\hspace{0 mm}
	\subfigure[][$\bV^d(4)=(0,0)$]{\includegraphics[width=4.2cm, height=4.2cm]{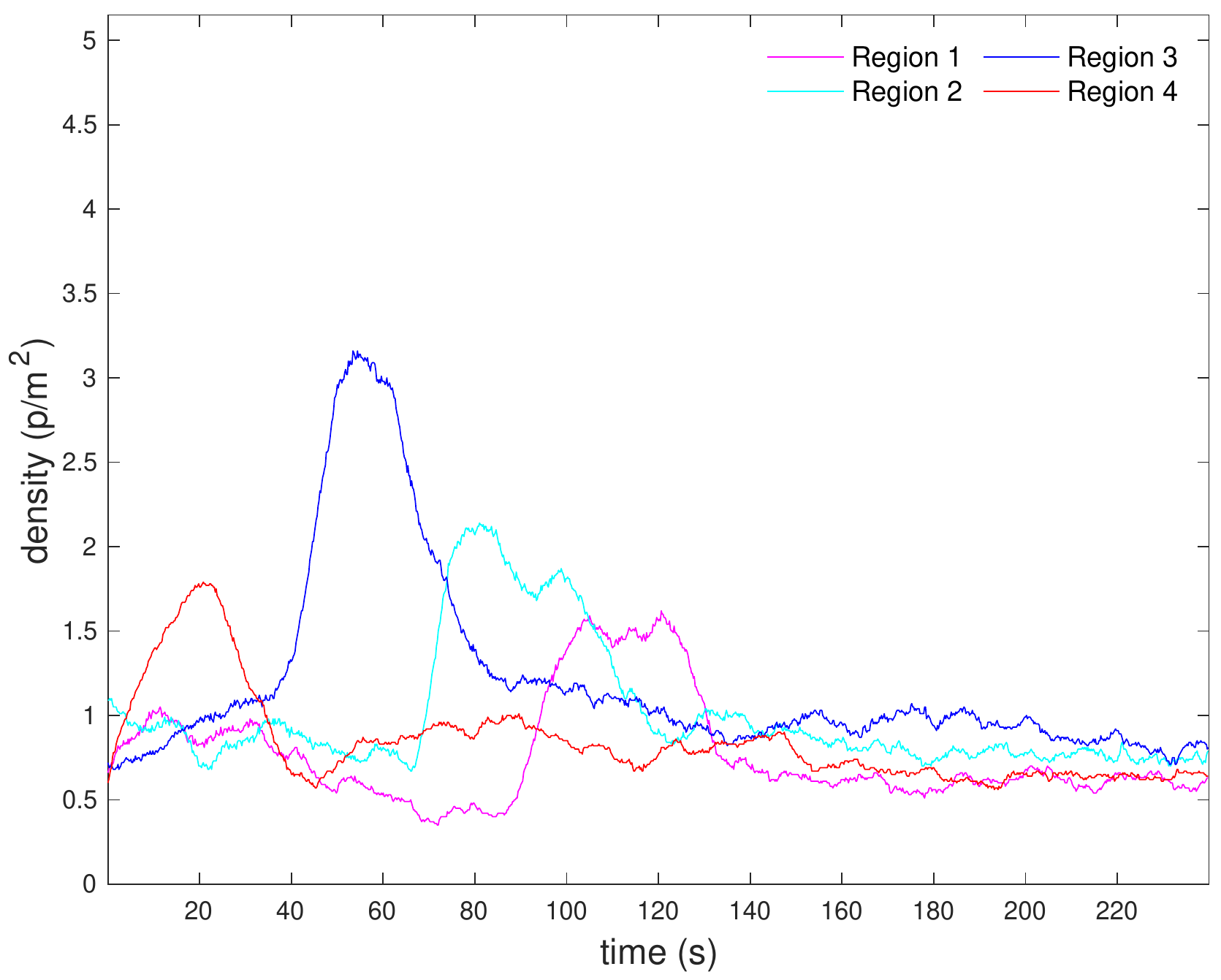}}
	\caption{Same simulation as in Fig.\ \ref{fig:densities_new}(c) with (a) $H$ decreased from 10 to 7 m, (b) $H$ increased from 10 to 13 m, (c) $D$ increased from 10 to 40 s, (d) $\delta t$ increased from 7 to 20 s, (e) $p$ increased from 25 to 75, (f) $\bV^d(4)$ decreased from $(0.5,0)$ to $(0,0)$ m/s.}
	\label{fig:densities_aggiunta}
\end{figure}

Now we consider the numerical setting of Fig.\ \ref{fig:densities_new}(c) as \emph{reference case} for investigating the role of the variable parameters of the model.
The results are shown in Fig.\ \ref{fig:densities_aggiunta}.

Figs.\ \ref{fig:densities_aggiunta}(a,b) show the result obtained modifying the height $H$ of the corridor. In both cases we maintain the initial density of people as 0.8 p/m$^2$. As expected, we observe the largest average densities in the case of the smallest corridor. The difference is remarkable in the regions farthest from the gate.

Fig.\ \ref{fig:densities_aggiunta}(c) shows the results obtained for an increased value of the doubt phase, namely $D=40$ s. 
We do not observe a relevant increase of the maximum values of the densities, but the peaks are more delayed.

Fig.\ \ref{fig:densities_aggiunta}(d) shows the results obtained for an increased value of the time needed to decide either to enter the doubt phase or to keep going to the gate, namely $\delta t=20$ s.
People are now more resistant to the idea of renouncing going to the gate and therefore necessarily spend more time in the corridor. Here the peak of density in region 4 is much greater than before and reaches the dangerous level of 5 p/m$^2$ \cite{helbing2012crowd}. 
This is mainly due to the longer time spent there by people in statuses 1--4 all together. 
Overall, the parameter $\delta t$ is the most effective in increasing the density.

Fig.\ \ref{fig:densities_aggiunta}(e) shows the effect of increasing $p$ from 25 to 75. In this case 75\% people decide to leave the corridor and only 25\% to stay. 
We see that at the beginning the densities in all regions are larger than those in the reference solution, because more people are trying to go back and form a big group which move very slowly. Conversely, at final time all regions are almost empty and have similar density.

Fig.\ \ref{fig:densities_aggiunta}(f) shows the effect of vanishing the desired velocity of people with status 4. Although those people are just waiting in the corridor, a rightward constant desired velocity is needed to reproduce the will to stay close to the gate and not to be overcome by newly incoming people in status 1. 
If $\bV^d(4)=(0,0)$ the dynamics change completely: people in status 1 continuously take the place of waiting people, and waiting people move back. 
At final time, the result is an almost constant density all along the corridor.

\subsection{Congestion control via Control \& Information Points}
\begin{figure}[h!]
	\centering 
	\subfigure[][Single CIP in region 3]{\includegraphics[width=4.2cm, height=4.2cm]{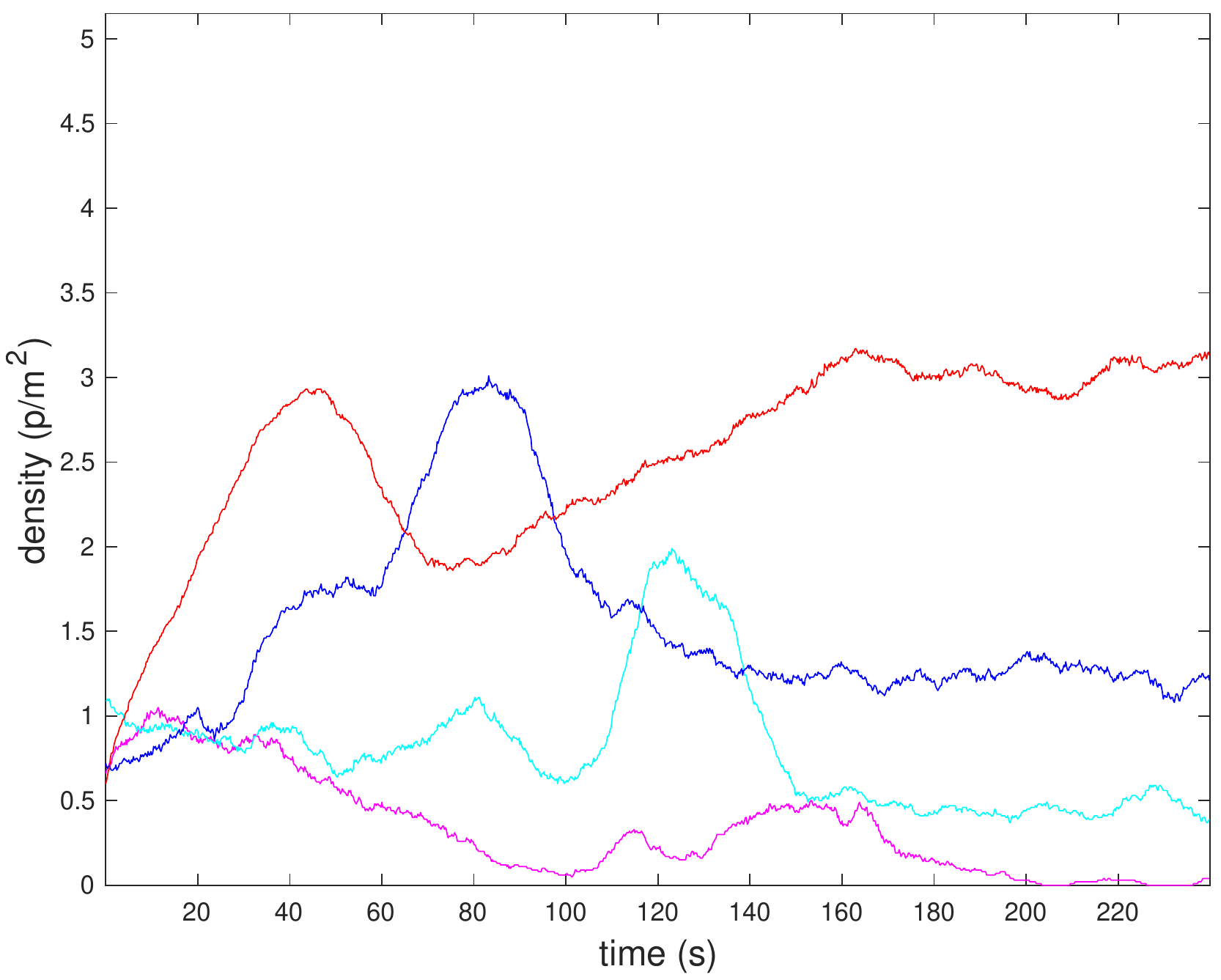}}
	\hspace{25 mm}
	\subfigure[][Two CIPs in regions 3 and 1]{\includegraphics[width=4.2cm, height=4.2cm]{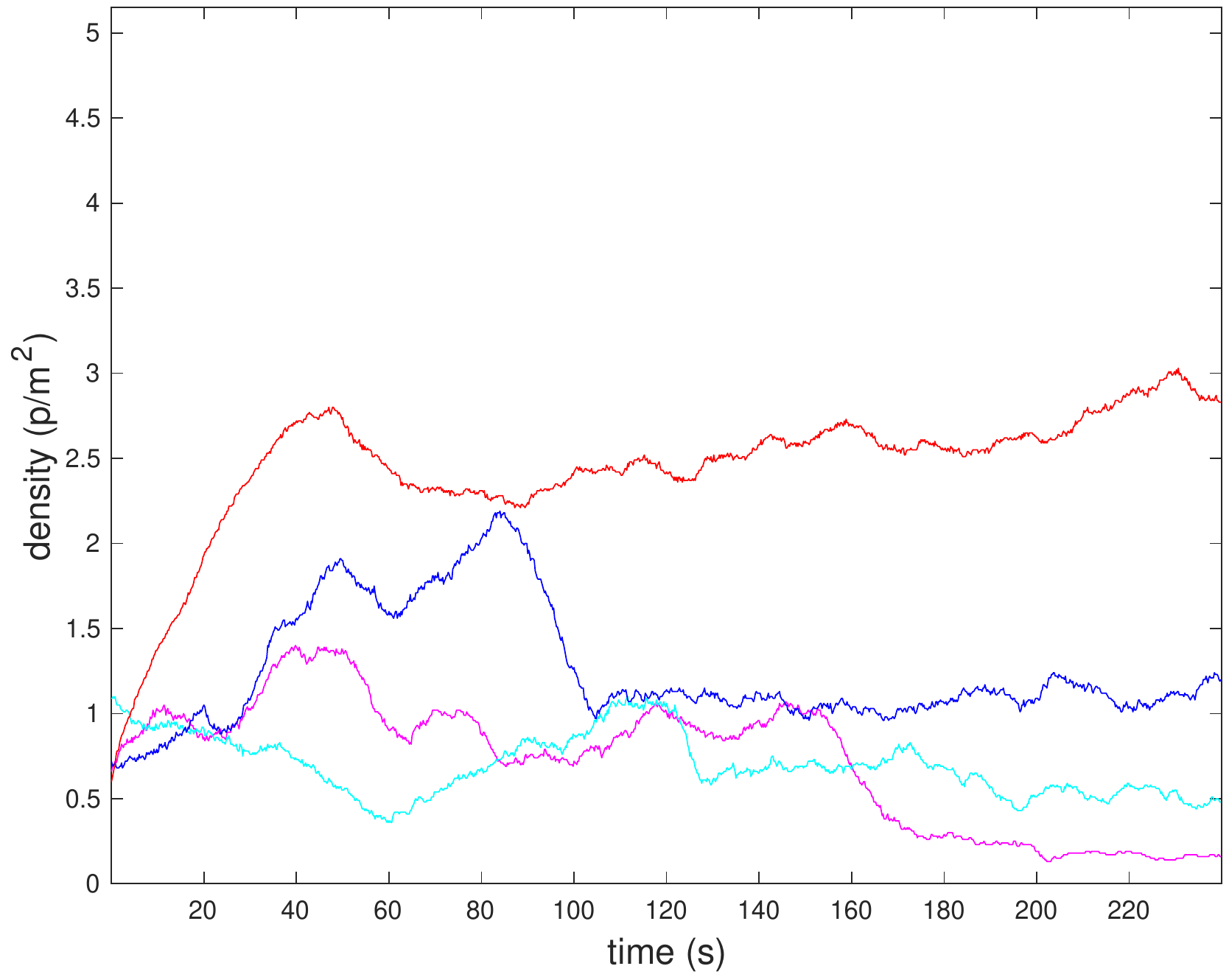}}
	\caption{Density evolution in the four regions of interest. 20 s after the gate closed, pedestrians are informed by (a) one CIP located in region 3, (b) two CIPs, simultaneously activated, located in regions 3 and 1.}
	\label{fig:signals}
\end{figure}

As final tests, we try to lower the congestion in the four regions adding some Control \& Information Points (CIPs) along the corridor. 
The idea is that people walking in the corridor are informed about the closure of the gate by some signals or stewards. 
In this way, people can enter the doubt phase much before (both in space and time) than they do without any knowledge of the status of gate. 
More precisely, we assume that anyone in status 1 moves immediately to status 2 as soon as it crosses a CIP.
The question arises when and where it is most convenient to locate CIPs.
Starting from the setting shown in Fig.\ \ref{fig:densities_aggiunta}(d), which is the most critical among the investigated ones, we have run a brute-force optimization procedure to test the effect of the presence of a single CIP; we have found that the best option is positioning the CIP in region 3, 20 s after the gate closed. 
The result is shown in Fig.\ \ref{fig:signals}(a). 

The effect of the CIP is clear since the maximal density decreases from $\sim$5 to $\sim$3 p/m$^2$, and the density evolution of region 3 becomes very similar to that of region 4.

Fig.\ \ref{fig:signals}(b) shows the result after adding a second CIP in region 1, again 20 s after the gate closed.
Maximal densities are not further decreased but densities are less fluctuating in all regions. This can be an advantage in terms of pedestrian safety.

\section{Conclusions and future perspectives}
We have conducted a numerical investigation in a particular situation characterized by a complex self-interaction of a crowd in a corridor. 
Although the model has a large number of parameters, we have seen that it is quite sensitive only to some of them. 
In particular, the parameter $N$ does not affect very much the maximal densities, as if the crowd was able to self-regulate. This ability comes from the dynamical way people change status, the decision being depending on the crowd itself.
The parameter $\delta t$, instead, affects the maximal density near the gate more than the size $N$ of the crowd itself. 
This conveys the idea that large crowds are not an issue \emph{per se}, but high densities actually arise whenever people keep moving towards a blocked crowd. This creates a compression which is then very hard to resolve. 

Finally, numerical simulations suggest that it is possible to prevent the formation of congestion by informing people along the corridor about the status of the gate. Even one information point is able to drastically reduce the maximal density since it prevents the encounter between the first people who have reached the closed gate and those who are arriving. 

\medskip

In the next future we will further investigate the role of the geometry of the domain in crowd dynamics, in particular testing different shapes of the corridor. 
We will also evaluate the impact of the visibility of the signals which inform people about the status of the gate.

Even more important, we will validate the results of the simulations comparing them with the real dynamics of people observed during mass events by means of video cameras.

Overall, we are strongly convinced that the adoption of predictive tools for congestion formation and technologies/automatisms for informing the crowd should be a further subject of study in event planning and crowd safety management.

\section*{Funding}
This work was carried out within the research project ``SMARTOUR: Intelligent Platform for Tourism'' (No. SCN\_00166) funded by the Ministry of University and Research with the Regional Development Fund of European Union (PON Research and Competitiveness 2007--2013). 

E.C. would also like to thank the Italian Ministry of Instruction, University and Research (MIUR) to support this research with funds coming from PRIN Project 2017 (No. 2017KKJP4X entitled ``Innovative numerical methods for evolutionary partial differential equations and applications'').

E.C.\ and M.M.\ are members of the INdAM Research group GNCS.

\section*{Authors' contribution}
G.A. proposed the research topic, suggested the case study, provided some real data for calibration, interpreted the results, and proofread the paper.
E.C. and M.M. developed the model, wrote the numerical code, performed numerical tests, interpreted the results, and wrote the paper.


\baselineskip=0.9\normalbaselineskip
\phantomsection\addcontentsline{toc}{section}{\numberline{}References}

\end{document}